\documentclass[11pt]{article}
\usepackage[margin=1in]{geometry}
\usepackage[T1]{fontenc}
\usepackage{lmodern,microtype}
\usepackage{amsmath,amssymb,amsthm,mathtools}
\usepackage{tikz}
\usepackage{float}
\usepackage[numbers,sort&compress]{natbib}
\usepackage[hidelinks]{hyperref}
\hypersetup{pdftitle={Exact chromatic numbers of small spheres}}
\usepackage{refcount}

\theoremstyle{plain}
\newtheorem{theorem}{Theorem}[section]
\newtheorem{lemma}[theorem]{Lemma}
\newtheorem{proposition}[theorem]{Proposition}
\newtheorem{corollary}[theorem]{Corollary}
\newtheorem{mainprime}{Theorem}

\theoremstyle{remark}
\newtheorem{remark}[theorem]{Remark}

\newcommand{\R}{\mathbb R}
\newcommand{\mc}{\mathcal M}
\newcommand{\ip}[2]{\langle #1,#2\rangle}

\usepackage{authblk}

\title{Exact chromatic numbers of small spheres}
\author[1]{Danila Cherkashin}
\author[2,3]{Vsevolod Voronov}

\affil[1]{\small Institute of Mathematics and Informatics, Bulgarian Academy of Sciences, Sofia}
\affil[2]{\small Laboratory of Advanced
Combinatorics and Network Applications,
Moscow Institute of Physics and Technology, Moscow, Russia}
\affil[3]{\small Caucasus Mathematical Center of Adyghe State University,
Maykop, Republic of Adygea, Russia}

\date{}
\begin{document}
\maketitle
\begin{abstract}
We prove that the unit-distance graph on the sphere $S^{n-1}(r)\subset\R^n$
has chromatic number at least $n+1$ for every $n \geq 3$ and $r>1/2$; for $r \leq r_\Delta(n)$ this bound is tight.
From a topological point of view this result means that for $d$ smaller than the diameter the closedness condition in the Hopf covering theorem for the sphere can be omitted.

Then for each such sphere we prove the existence of a finite $(n+1)$-chromatic graph that admits
a unit-distance realization on the image of every sufficiently small continuous perturbation of the sphere. 
\end{abstract}

\section{Introduction}
Let $S^{n-1}(r) = \{x\in\R^n:\|x\|=r\}$ be the sphere of radius $r$ centered at the origin. 
Define the unit-distance graph on a set $M\subset\R^n$ by $(M,E)$, where $E$ is the set of pairs from $M$ at Euclidean
distance one. 
Let $\chi(M)$ denote the chromatic number of this graph, which is the smallest number of colors one needs to color points of $M$ such that no monochromatic pair forms an edge.

Lov\'asz~\cite{Lovasz1983} proved that $\chi(S^{n-1}(r)) \geq n$ for
every $n \geq 2$ and every $r>1/2$ and obtained the lower bound $n+1$ along a sequence of
radii tending to $1/2$. He asked whether $n$ colors ever suffice when
$n>2$. For $n=3$, Simmons~\cite{Simmons1976} proved that four colors are
necessary for $r\geq 1/\sqrt3$. In~\cite{CherkashinVoronov2024} the authors extended this result to every $r>1/2$.
In this paper we answer Lov\'asz's question in all dimensions $n\geq 3$.

\begin{theorem}\label{thm:main}
For every integer $n\geq 3$ and every $r>1/2$,
\[
\chi\bigl(S^{n-1}(r)\bigr)\geq n+1.
\]
\end{theorem}

Theorem~\ref{thm:main} is tight for some open interval of $r$. 
Let $\Delta$ be a regular $n$-simplex inscribed in
$S^{n-1}(r)$. The radial projections of its facets from the
origin form a cover of the sphere by $n+1$ congruent closed cells.
Assign distinct colors to the cells and give each point the
color of one cell containing it, choosing arbitrarily when
several cells contain the point.
This coloring is proper for $r < r_\Delta(n)$, see~\cite[Lemma~1]{Santalo1946} and~\cite{Horvath2021}, where
\[
r_\Delta(n) :=
\begin{cases}
\displaystyle \frac12\sqrt{\frac{n+2}{n+1}}, & n \text{ even},\\[8pt]
\displaystyle \frac{1}{\sqrt{\,2\left(1+\sqrt{\frac{n-1}{n+3}}\right)}}, & n \text{ odd}.
\end{cases}
\]

The de Bruijn--Erd\H{o}s theorem~\cite{deBruijnErdos1951} also gives the upper bound at $r = r_\Delta(n)$. 
Then we have the following
\begin{corollary} \label{cor:boundary} 
For every $n\geq 3$ and every $r$ with $1/2<r\leq r_\Delta(n)$,
\[
\chi\bigl(S^{n-1}(r)\bigr)=n+1.
\]
\end{corollary}

The interval of radii for which $\chi(S^{n-1}(r))=n+1$ may extend beyond $r_\Delta(n)$. In particular, for $n\geq 4$ it is not known for which $r$ the sphere $S^{n-1}(r)$ can be partitioned into $n+1$ parts of diameter less than one~\cite{kalai2015some}.

The endpoint $r_\Delta(n)$ tends to $1/2$ as $n\to\infty$.
Raigorodskii~\cite{Raigorodskii2012} proved that for every fixed
$r>1/2$, the chromatic number $\chi(S^{n-1}(r))$ grows exponentially
with $n$. Raigorodskii and Sagdeev~\cite{RaigorodskiiSagdeev2026} showed that there is a constant $C>0$ such that,
for all sufficiently large $n$, the inequality $r\ge 1/2 + C \sqrt{\ln n / n} $ implies $\chi(S^{n-1}(r))>n+1$.

Borsuk~\cite{Borsuk1933} proved that if $m+1$ closed sets cover $S^m$,
one contains an antipodal pair. Hopf~\cite{Hopf1944} proved that for
each prescribed chord length $0 < d \leq 2$, some member of the cover contains
a pair at distance $d$. Bogatyi~\cite[\S1]{Bogatyi1986} discusses this
covering theorem and its extensions to other configurations.
Theorem~\ref{thm:main} is equivalent to the following version for arbitrary covers.

\begin{mainprime}[Hopf covering theorem without closedness]\label{cor:hopf}
Let $n\ge 3$ and suppose 
\[
S^{n-1}(1)=A_1\cup\cdots\cup A_n.
\]
For every $0<d<2$, some $A_i$ contains two points at distance $d$.
\end{mainprime}

Rescaling the sphere gives a convenient way to state the stability result.
Put $S = S^{n-1}(1)$.
For $x,y\in S$, the points $rx$ and $ry$ are at distance one exactly when
$\ip{x}{y}=a$, where $a=1-1/(2r^2)$. As $r$ ranges over $(1/2,\infty)$,
$a$ ranges over $(-1,1)$. For a continuous symmetric kernel function
$K:S\times S\to\R$, let $G(K)$ be the graph with vertex set $S$ whose edges join distinct points $x,y$ when $K(x,y)=0$.
Write $K_a(x,y)=\ip{x}{y}-a$ and $G_a=G(K_a)$.

\begin{theorem}\label{thm:transfer}
Let $n\geq 3$. For every $a \in (-1,1)$ there is $\varepsilon>0$ such that every continuous symmetric $K$ with $\|K-K_a\|_\infty<\varepsilon$ satisfies $\chi(G(K)) \geq n+1$.
\end{theorem}

\paragraph{Structure of the paper.} 
Sections~2 and~3 are devoted to the proof of Theorem~\ref{thm:main}: Section~2 shows that there are no proper $n$-colorings with dense color classes. 
For $1/2<r<1/\sqrt2$, Section~3 provides the complementary statement.
The remaining radii are handled by induction on the dimension.
The main tool in the proof is the analysis of perturbed embeddings of two sequences of graphs: odd cycles with attached pendant vertices and generalized Mycielski graphs without the apex. Section~\ref{sec:transfer} proves stability of the chromatic lower bound (Theorem~\ref{thm:transfer})
by taking limits of color classes. Section~\ref{sec:regularwitness} constructs a single finite graph (Theorem~\ref{thm:stableintro}) that witnesses the bound for every
sufficiently small perturbation at a given $a$, strengthening Theorem~\ref{thm:transfer}.

\section{Excluding dense colorings via pinned odd cycles} \label{sec:tools}

For a smooth curve $y:[0,t_1]\to S$, we consider its time parametrization $y(t)$ and denote its derivative by
$\dot y=dy/dt$. We write a family indexed by $\mathbb Z_m$ as $y=(y_0,\ldots,y_{m-1})$, with subscripts read modulo $m$.
For a smooth manifold $M$, let $T_xM$ denote its tangent space at $x$. The tangent bundle of $M$ is
\[
TM=\{(x,v):x\in M,\ v\in T_xM\},
\]
with projection $(x,v)\mapsto x$.
For a smooth map $F:M\to N$, its differential at $x$ is the linear map
\[
DF(x):T_xM\longrightarrow T_{F(x)}N
\]
characterized by $DF(x)\gamma'(0)=(F\circ\gamma)'(0)$ for every smooth
curve $\gamma$ in $M$ with $\gamma(0)=x$. These maps define
$DF:TM\to TN$ by $(x,v)\mapsto(F(x),DF(x)v)$.
When $N=\mathbb R^q$, we identify $T_{F(x)}\mathbb R^q$ with $\mathbb R^q$.
We will need the following technical statement about zeros of perturbed maps.

\begin{lemma}\label{lem:zero}
Let $M$ be a smooth manifold, let $U\subset M$ be open, and let
$F:U\to\R^q$ be smooth. Suppose $u\in U$, $F(u)=0$, and
$DF(u):T_uM\to\R^q$ is surjective. For every open neighborhood
$V\subset U$ of $u$, there exists $\delta>0$ such that every continuous
map $\widetilde F:V\to\R^q$ satisfying
\[
 \sup_{v\in V}\|\widetilde F(v)-F(v)\|<\delta
\]
has a zero in $V$.

In particular, if continuous maps $F_j:U\to\R^q$ converge uniformly
to $F$ on a neighborhood of $u$, then for some $j_0$ there are points
$u_j\in U$ with $F_j(u_j)=0$ for all $j\geq j_0$ and $u_j\to u$.
\end{lemma}
\begin{proof}
Restrict $F$ to a $q$-dimensional disk through $u$ with tangent space
complementary to $\ker DF(u)$. The inverse function theorem gives a
local inverse $g$ with $g(0)=u$. Choose $\delta>0$ so that $g$ is
defined on the closed ball $B$ of radius $\delta$ centered at $0$ and
$g(B)\subset V$. Since $F(g(z))=z$ on $B$, the map
$z\mapsto z-\widetilde F(g(z))$ sends $B$ into itself.
By Brouwer's fixed point theorem it has a fixed point $z$.
Then $\widetilde F(g(z))=0$ and $\|z\|=\|\widetilde F(g(z))-F(g(z))\|$.
For the sequence $F_j$, choose $V$ inside a neighborhood of uniform
convergence. For large $j$, the uniform error on $V$ is less than
$\delta$, and the corresponding fixed points $z_j$ tend to $0$.
Thus $u_j=g(z_j)\to u$.
\end{proof}

\begin{figure}[H]
\centering
\begin{tikzpicture}[font=\small]
 \pgfmathsetmacro{\cycleradius}{1.9}
 \pgfmathsetmacro{\edgelength}{2*\cycleradius*sin(72)}
 \foreach \s in {0,...,4} {
 \coordinate (Y\s) at ({90+144*\s}:\cycleradius);
 }
 \foreach \s in {0,...,4} {
 \pgfmathsetmacro{\vertexangle}{90+144*\s}
 \foreach \j/\offset in {1/-32,2/32} {
 \draw[black!55] (Y\s) -- ++({\vertexangle+180+\offset}:\edgelength)
 coordinate (P\s-\j);
 }
 }
 \draw[thick] (Y0)--(Y1)--(Y2)--(Y3)--(Y4)--cycle;
 \foreach \s in {0,...,4} {
 \foreach \j in {1,2} {
 \draw[fill=white] (P\s-\j) circle (2pt);
 }
 \fill (Y\s) circle (2pt);
 \node at ({90+144*\s}:{\cycleradius+.28}) {$y_{\s}$};
 }
 \node[left] at (P0-1) {$x_{0,1}$};
 \node[right] at (P0-2) {$x_{0,2}$};
\end{tikzpicture}
\caption{Schematic of a five-cycle with two pins (open circles) at
each vertex.}
\label{fig:five-cycle-pins}
\end{figure}
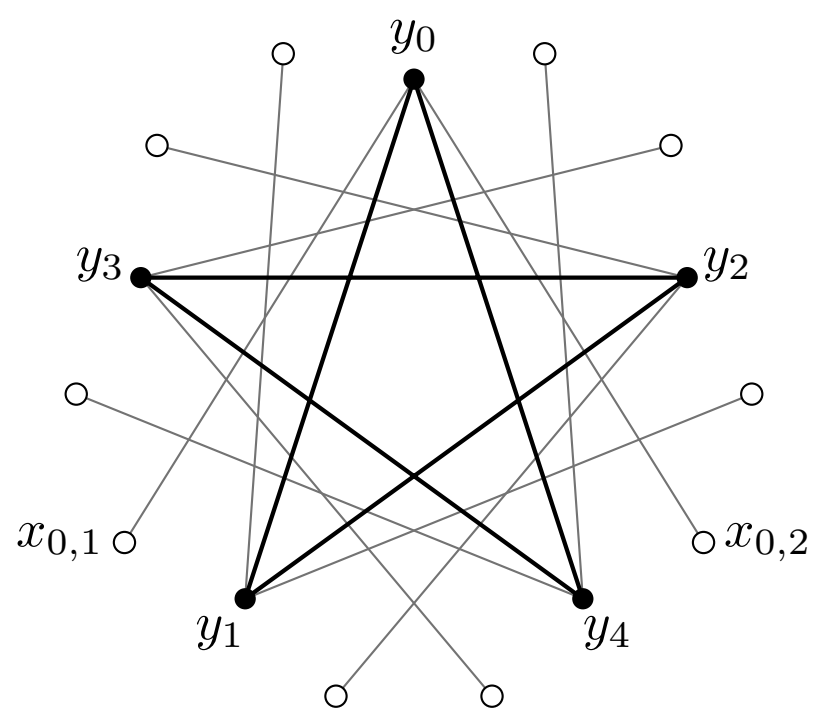

We attach $n-2$ vertices, called \emph{pins}, to each vertex of an odd
cycle (Figure~\ref{fig:five-cycle-pins}). The case of one pin on $S^2$
was considered in~\cite{CherkashinVoronov2024}.

\begin{lemma}[A pinned odd cycle]\label{lem:pinnedcycle}
For $n\geq3$ and $-1<a<1$, there exist an odd integer $m\geq3$,
cycle positions $y_s\in S$, and pin positions $x_{s,j}\in S$ satisfying
\[
 \ip{y_s}{y_{s+1}}=a,\qquad \ip{y_s}{x_{s,j}}=a \qquad \text{for}
 \quad s\in \mathbb Z_m,\ 1\leq j\leq n-2.
\]
With the pins fixed, the map $S^{\times m}\to\R^{m(n-1)}$ given by
the left-hand sides has an invertible differential at
$(y_0,\ldots,y_{m-1})$.
In particular, there are open neighborhoods $U_{s,j}$ of $x_{s,j}$
such that the cycle positions solving these equations depend smoothly
on the pin positions in $\prod_{s,j}U_{s,j}$.
\end{lemma}
\begin{proof}
Choose an odd $m\geq3$ with $b=\cos(\pi/m)>-a$. In a copy of $\R^3$
inside $\R^n$, put
\[
 \alpha=\pi-\pi/m,\qquad h=\sqrt{(b+a)/(1+b)},\qquad
 y_s=(\sqrt{1-h^2}\cos(s\alpha),\sqrt{1-h^2}\sin(s\alpha),h).
\]
Then $0<h<1$ and $\ip{y_s}{y_{s+1}}=h^2-(1-h^2)b=a$.
The horizontal projections are distinct since $(m-1)/2$ and $m$ are coprime.
The projections lie on a circle, so any three are noncollinear.
Since all the points $y_s$ have the same nonzero height, any three of
their position vectors are linearly independent.

Let $v_s$ be the normalized projection of $y_{s+1}$ onto
the orthogonal complement $\operatorname{span}(y_{s-1},y_s)^\perp$. Then $v_s\perp y_{s-1},y_s$
and $c_s:=\ip{v_s}{y_{s+1}}>0$.
For an orthonormal basis $(w_{s,j})_{j=1}^{n-2}$ of
$\operatorname{span}(y_s,v_s)^\perp$, take the pins
\[
 x_{s,j}=a y_s+\sqrt{1-a^2}\,w_{s,j}.
\]
These are unit vectors with $\ip{y_s}{x_{s,j}}=a$.

Let $(\xi_s)_s$, with $\xi_s\in T_{y_s}S$, lie in the kernel of the
differential with the pins fixed. The pin equations give
$\ip{\xi_s}{x_{s,j}}=0$ for every $s,j$, so $\xi_s=\lambda_sv_s$.
The cycle equations give
\[
 0=\lambda_s\ip{v_s}{y_{s+1}}+\lambda_{s+1}\ip{y_s}{v_{s+1}}
 =\lambda_sc_s.
\]
Since $c_s>0$, the kernel is zero. There are $m(n-1)$ equations and as
many tangent coordinates, so the differential is invertible. The
implicit function theorem gives the final assertion.
\end{proof}

\begin{corollary}\label{lem:cycle}
For $n\geq3$ and $-1<a<1$, $G_a$ has no proper $n$-coloring in which
the first $n-2$ color classes are dense.
\end{corollary}
\begin{proof}
In such a coloring, choose a point of color $j$ in each $U_{s,j}$
from Lemma~\ref{lem:pinnedcycle}. Each vertex of the resulting cycle
is adjacent to a pin of every color $1,\ldots,n-2$. The remaining two
colors cannot properly color the odd cycle.
\end{proof}

\section{Spherical embeddings of Mycielski graphs}

For small radii, the geometric argument from~\cite{CherkashinVoronov2024} shows that every
color class of a hypothetical proper three-coloring of the unit-distance graph on $S^2(r)$ would be dense.
In higher dimensions, we use realizations of generalized Mycielski
graphs at a dense set of edge inner products in one lower dimension.
This lets us place the last layer of a generalized Mycielski graph 
in any prescribed nonempty open set. 

For $1/2<r<1/\sqrt2$, the parameter $a = 1-1/(2r^2)$ is negative. 
A realization assigns a point to each vertex so that every edge has
the prescribed inner product. Nonadjacent vertices may share a point.
Any proper coloring of the image graph gives a proper coloring of the
original graph.

The generalized Mycielski graph $M_k(H)$ has $k\geq1$ copies of $V(H)$,
called layers $0,\ldots,k-1$, and an apex. Put a copy of $H$ in layer
$0$. For each $uv\in E(H)$ and $0\leq j<k-1$, join $(u,j)$ to
$(v,j+1)$ and $(v,j)$ to $(u,j+1)$. Join the apex to all vertices in
layer $k-1$. The case $k=2$ is Mycielski's original construction
\cite{Mycielski1955}. For the general construction see
Sachs and Stiebitz~\cite{SachsStiebitz1989}. For $m\geq2$, let $\mc_m$
consist of the graphs obtained from $K_2$ by $m-2$ such steps, choosing
$k$ independently at each step. By Stiebitz's theorem
\cite{Stiebitz1985,MullerStehlik2019},
\begin{equation}\label{eq:stiebitz}
 H\in\mc_m\quad\Longrightarrow\quad\chi(H)\geq m.
\end{equation}

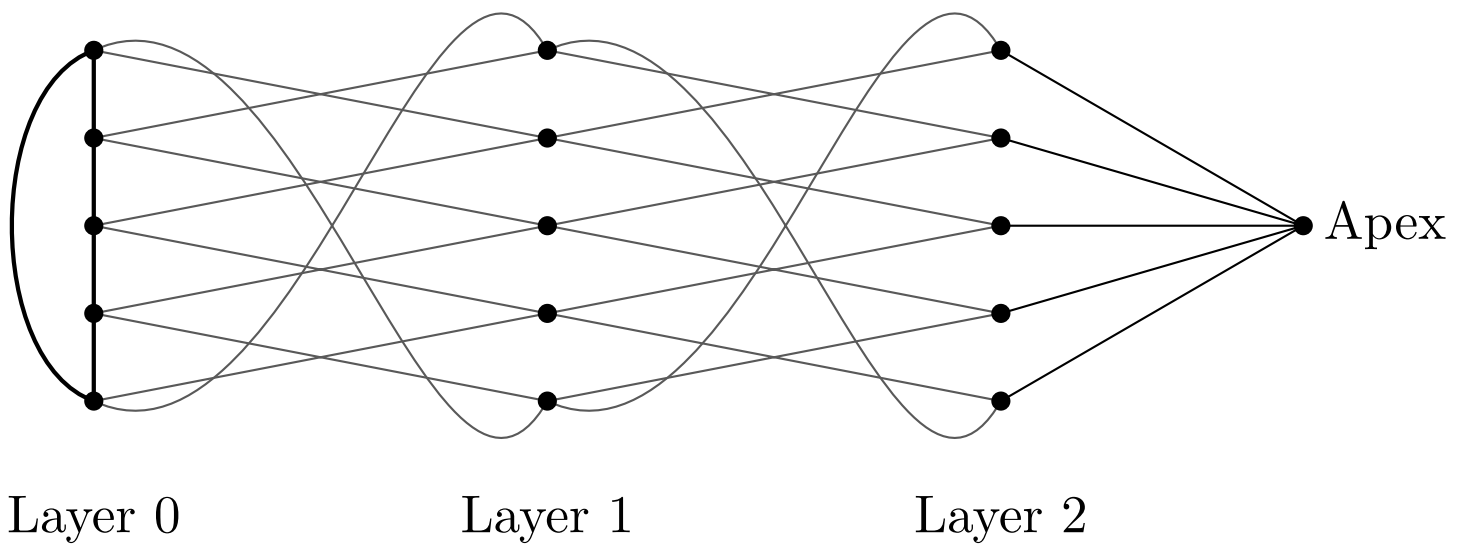
\begin{figure}[H]
\centering
\begin{tikzpicture}[font=\small]
 \foreach \layer in {0,1,2} {
 \foreach \vertex in {0,...,4} {
 \coordinate (L\layer-\vertex) at ({3*\layer},{.58*(2-\vertex)});
 }
 }
 \coordinate (apex) at (8,0);

 \foreach \layer/\nextlayer in {0/1,1/2} {
 \foreach \vertex/\nextvertex in {0/1,1/2,2/3,3/4} {
 \draw[black!65] (L\layer-\vertex)--(L\nextlayer-\nextvertex)
 (L\layer-\nextvertex)--(L\nextlayer-\vertex);
 }
 \draw[black!65] (L\layer-0) to[out=25,in=-120] (L\nextlayer-4)
 (L\layer-4) to[out=-25,in=120] (L\nextlayer-0);
 }
 \draw[thick] (L0-0)--(L0-1)--(L0-2)--(L0-3)--(L0-4)
 to[out=160,in=200,looseness=.85] (L0-0);
 \foreach \vertex in {0,...,4} {
 \draw (L2-\vertex)--(apex);
 }
 \foreach \layer in {0,1,2} {
 \foreach \vertex in {0,...,4} {
 \fill (L\layer-\vertex) circle (1.8pt);
 }
 \node[below] at ({3*\layer},-1.65) {Layer \layer};
 }
 \fill (apex) circle (1.8pt) node[right] {Apex};
\end{tikzpicture}
\caption{The generalized Mycielski graph $M_3(C_5)$.}
\label{fig:mycielski-five-cycle}
\end{figure}

In any proper three-coloring of the graph in
Figure~\ref{fig:mycielski-five-cycle} with its apex deleted, every
layer uses all three colors. If layer $j$ missed a color, keeping
layers $0,\ldots,j$ and giving the new apex that color would yield
a proper three-coloring of $M_{j+1}(C_5)$. This contradicts
\eqref{eq:stiebitz}, since $C_5=M_2(K_2)$.

This layered construction occurs implicitly in Lov\'asz's proof
\cite[Theorem~1, Claims~1 and~3]{Lovasz1983}. M\"uller and Stehl\'ik
\cite[Lemma~7]{MullerStehlik2019} used alternating latitudes to
place adjacent vertices nearly antipodally. We choose layer angles
that give each edge inner product $a$.

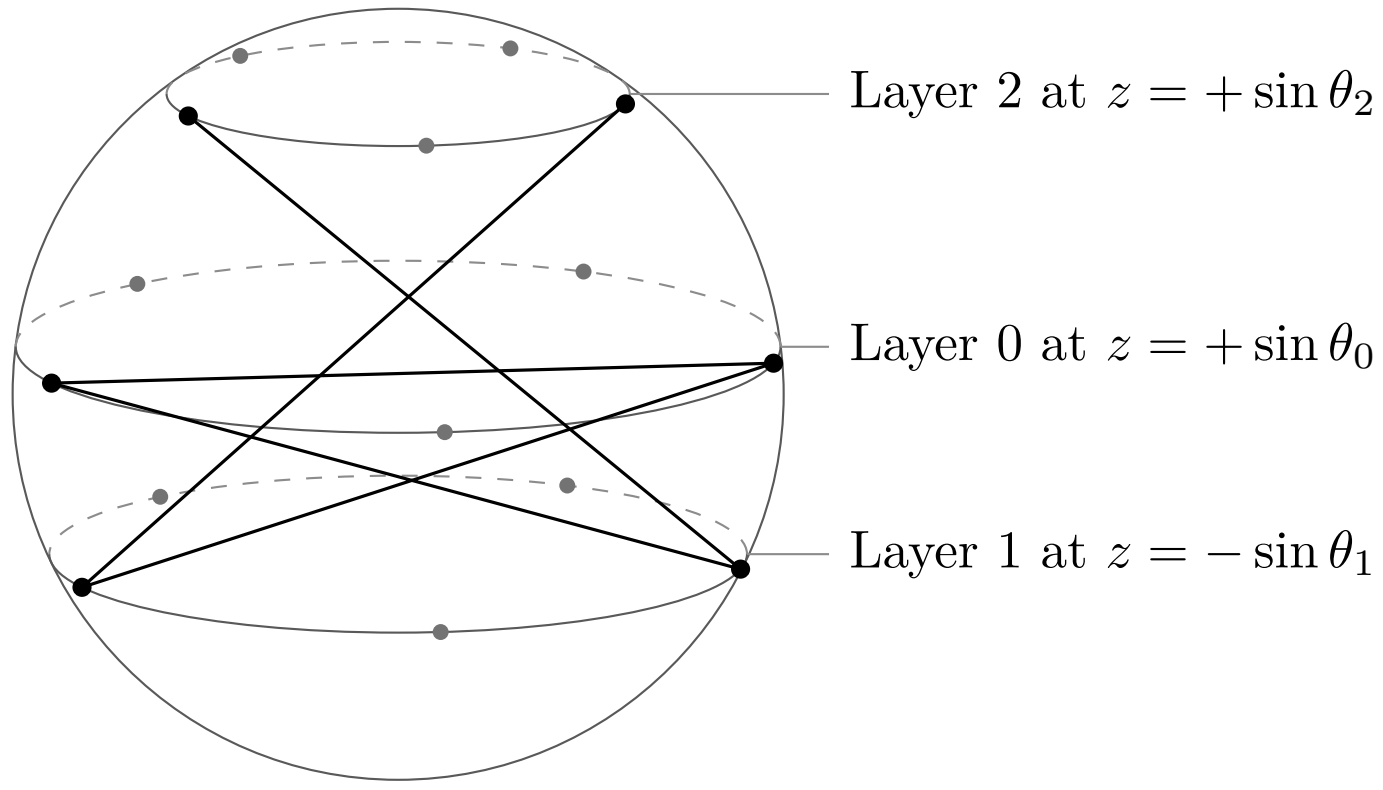
\begin{figure}[H]
\centering
\begin{tikzpicture}[font=\small]
 \pgfmathsetmacro{\sphereradius}{2.55}
 \pgfmathsetmacro{\viewangle}{13}
 \pgfmathsetmacro{\bvalue}{cos(36)}
 \pgfmathsetmacro{\pvalue}{.78}
 \pgfmathsetmacro{\anglezero}{asin(sqrt((\bvalue-\pvalue)/(1+\bvalue)))}
 % Take the root of the cross-layer equation after the sinusoid's maximum.
 \pgfmathdeclarefunction{nextlatitude}{1}{%
 \pgfmathparse{atan2(sin(#1),\bvalue*cos(#1))
 + acos(\pvalue/sqrt((\bvalue*cos(#1))^2+sin(#1)^2))}%
 }
 \pgfmathsetmacro{\angleone}{nextlatitude(\anglezero)}
 \pgfmathsetmacro{\angletwo}{nextlatitude(\angleone)}

 \draw[black!65] (0,0) circle (\sphereradius);
 % Project (x,y,z) to (x,y sin(view angle) + z cos(view angle)).
 \foreach \layer/\latitude/\sign in {0/\anglezero/+,1/-\angleone/-,2/\angletwo/+} {
 \pgfmathsetmacro{\height}{\sphereradius*sin(\latitude)*cos(\viewangle)}
 \pgfmathsetmacro{\majoraxis}{\sphereradius*cos(\latitude)}
 \pgfmathsetmacro{\minoraxis}{\majoraxis*sin(\viewangle)}
 \draw[black!45,dashed] (\majoraxis,\height)
 arc[start angle=0,end angle=180,x radius=\majoraxis,y radius=\minoraxis];
 \draw[black!65] (-\majoraxis,\height)
 arc[start angle=180,end angle=360,x radius=\majoraxis,y radius=\minoraxis];
 \draw[black!45] (\majoraxis,\height) -- (2.85,\height);
 \node[right] at (2.85,\height) {Layer $\layer$ at $z=\sign\sin\theta_{\layer}$};
 \foreach \vertex in {0,...,4} {
 \pgfmathsetmacro{\azimuth}{205+144*\vertex}
 \coordinate (X\layer-\vertex) at
 ({\majoraxis*cos(\azimuth)},{\height+\minoraxis*sin(\azimuth)});
 }
 }

 \draw[semithick] (X0-0) -- (X0-1);
 \foreach \layer/\nextlayer in {0/1,1/2} {
 \draw[semithick] (X\layer-0) -- (X\nextlayer-1)
 (X\layer-1) -- (X\nextlayer-0);
 }
 \foreach \layer in {0,1,2} {
 \foreach \vertex in {2,3,4} {
 \fill[black!55] (X\layer-\vertex) circle (1.5pt);
 }
 \fill (X\layer-0) circle (1.8pt);
 \fill (X\layer-1) circle (1.8pt);
 }
% \node[above left,xshift=-3pt] at (X0-0) {$X_0(u)$};
 % \node[below right] at (X0-1) {$X_0(v)$};
\end{tikzpicture}
\caption{The latitudes of $M_3(C_5)$ without the apex.}
\label{fig:mycielski-latitudes}
\end{figure}

Suppose $H$ is realized on $S^{n-2}(1)$ at edge inner product $-b$,
where $0<-a<b<1$. For a vertex at $u$, place its copy in layer $j$ at
\begin{equation}\label{eq:layer}
 X_j(u)=(\cos\theta_j\,u,\,(-1)^j\sin\theta_j)\in S.
\end{equation}
The edges in layer $0$ and those between layers $j$ and $j+1$
have inner product $a$ when
\begin{align}
 b\cos^2\theta_0-\sin^2\theta_0 &= -a,\label{eq:start}\\
 b\cos\theta_j\cos\theta_{j+1}+\sin\theta_j\sin\theta_{j+1} &= -a
 \qquad(j\geq0).
 \label{eq:cross}
\end{align}
Figure~\ref{fig:mycielski-latitudes} illustrates this placement for $H=C_5$, $k=3$ and $n=3$.

If $\theta_k=\pi/2$ for some $k\geq1$, layer $k$ collapses to the pole
$(-1)^k e_n$, where $e_n$ is the last coordinate unit vector.
Then \eqref{eq:cross} gives $\sin\theta_{k-1}= -a$, so the pole has
inner product $a$ with every point of layer $k-1$. Layers
$0,\ldots,k-1$ and this apex realize $M_k(H)$.

\begin{lemma}\label{lem:shoot}
Equations \eqref{eq:start}--\eqref{eq:cross} have solutions
$\theta_j(b,-a)$ continuous on $0<-a<b<1$. For fixed $(b,-a)$,
$\theta_j\to+\infty$ as $j\to\infty$.
For fixed $j$, $\theta_j\to0$ as $(b,-a)\to(\beta,\beta)$,
$0<\beta<1$.
\end{lemma}
\begin{proof}
Take $\theta_0=\arcsin\sqrt{(b+a)/(1+b)}$. Given $x=\theta_j$, choose
$y=\theta_{j+1}$ in $(x,x+\pi)$ with $D_b(x,y)=-a$, where
\[
 D_b(x,y)=b\cos x\cos y+\sin x\sin y.
\]
As a function of $y$, $D_b(x,y)$ is a sinusoid with amplitude
$\sqrt{b^2\cos^2x+\sin^2x}\in[b,1]$. Its values at $x$ and $x+\pi$
straddle $-a$, since
$D_b(x,x)\geq b>-a$ and $D_b(x,x+\pi)=-D_b(x,x)$.
The crossings are simple, so their number in the interval is odd.
Three successive crossings span $2\pi$, so there is exactly one.
It depends continuously on $(b,-a,x)$ by the implicit function theorem.
Since $|\partial_y D_b|\leq1$, the mean value theorem gives
$y-x\geq D_b(x,x)+a\geq b+a$, so the sequence diverges.

As $(b,-a)\to(\beta,\beta)$, the formula gives $\theta_0\to0$.
If $\theta_j\to0$, then every limit point $y$ of $\theta_{j+1}$ lies
in $[0,\pi]$ and satisfies $\beta\cos y=\beta$ by \eqref{eq:cross}.
Thus $\theta_{j+1}\to0$ as well.
\end{proof}

\begin{proposition}\label{prop:dense}
For $n\geq2$, the set $D_n$ of parameters $a\in(-1,0)$ for which some
$H\in\mc_{n+1}$ has a realization on $S^{n-1}(1)$ at edge inner
product $a$ is dense in $(-1,0)$.
\end{proposition}
\begin{proof}
Fractions $h/q$ in lowest terms with $q$ odd are dense in $(1/4,1/2)$,
since one can approximate with arbitrarily large odd denominators and
then reduce. Stepping by $2\pi h/q$ on the circle realizes the odd
cycle $C_q=M_{(q-1)/2}(K_2)$ at $a=\cos(2\pi h/q)$. Thus $D_2$ is dense.

For $n\geq3$ and $0<u<v<1$, choose $b\in(u,v)$ such that $-b\in D_{n-1}$ by induction.
Let $H\in\mc_n$ be the corresponding graph. By Lemma~\ref{lem:shoot},
choose $k\geq1$ with $\theta_k(b,u)>\pi/2$. Since
$\theta_k(b,p)\to0$ as $p \to b^-$, there is a $p\in(u,b)$ with
$\theta_k(b,p)=\pi/2$. At this $p$, the pole gives a realization of
$M_k(H)\in\mc_{n+1}$, so $-p\in D_n\cap(-v,-u)$.
\end{proof}

The next lemma shows that every proper $n$-coloring of $G_a$ with $a < 0$ has dense color classes.

\begin{lemma}\label{lem:cap}
Let $n\geq3$, $a \in (-1,0)$, and let $U\subset S$ be nonempty and open.
There is a finite subgraph $F\subset G_{a}$ such that every proper $n$-coloring of $F$ uses
all $n$ colors on $V(F) \cap U$.
\end{lemma}
\begin{proof}
Choose $b_1\in(-a,1)$ and take $k\geq1$ with
$\theta_k(b_1,-a)>\pi/2$. Since $\theta_k(b,-a)\to0$ as $b \to -a^+$,
there is a $b_0\in(-a,b_1)$ with $\theta_k(b_0,-a)=\pi/2$.
By \eqref{eq:layer}, $X_k(u)$ tends uniformly in $u$ to $(-1)^k e_n$
as $b\to b_0$. Fix a rotation taking this pole into $U$.
By Proposition~\ref{prop:dense}, choose $b$ with $-b\in D_{n-1}$ so close to $b_0$ that the rotated
layer $k$ lies in $U$.

For the corresponding $H\in\mc_n$, equations
\eqref{eq:start}--\eqref{eq:cross} realize layers $0,\ldots,k$ of
$M_{k+1}(H)$. Let $F$ consist of the rotated images of their vertices
and edges, and let $A$ be the image of layer $k$. A proper $n$-coloring
of $F$ gives a proper coloring of the layers. Any color missing from
$A$ could be assigned to the apex, contradicting \eqref{eq:stiebitz}
for $M_{k+1}(H)\in\mc_{n+1}$.
\end{proof}

We are now ready to prove the main result.

\begin{proof}[Proof of Theorem~\ref{thm:main}]
First suppose that $1/2<r<1/\sqrt2$, so that $-1<a<0$.
In a proper $n$-coloring of $G_{a}$, Lemma~\ref{lem:cap} would force every color
class to meet every nonempty open set. Corollary~\ref{lem:cycle}
excludes such a coloring.

For $r\geq1/\sqrt2$, induct on $n$, starting with $n=3$
\cite{CherkashinVoronov2024}. The unit-distance neighbors of
$x\in S^{n-1}(r)$ form the section $\ip{x}{y}=r^2-1/2$, an
$(n-2)$-sphere of radius
\[
 s=\sqrt{r^2-(r^2-1/2)^2/r^2}
 =\sqrt{1-1/(4r^2)}\geq1/\sqrt2.
\]
For $n\geq4$, this sphere requires at least $n$ colors by induction,
and the color of $x$ is different from all of them.
\end{proof}

\begin{proof}[Proof of Corollary~\ref{cor:boundary}]
By Theorem~\ref{thm:main} and the coloring from the Introduction, it remains to show that
$\chi(S^{n-1}(r))\le n+1$ for $r=r_\Delta(n)$. In this case every closed cell has
diameter one, and two interior points of the same cell are at distance less than
one. Let $F\subset S^{n-1}(r)$ be finite. The
union of the cell boundaries is closed and nowhere dense, so some rotation $R$
maps every point of $F$ into the interior of a cell. Give each $p\in F$ the color
of the cell containing $R(p)$. If $\|p-q\|=1$, then $\|R(p)-R(q)\|=1$, so $R(p)$
and $R(q)$ lie in different cells. Thus every finite subgraph is
$(n+1)$-colorable, and the de~Bruijn--Erd\H{o}s theorem~\cite{deBruijnErdos1951} gives the claim.
\end{proof}
%For any finite subset of points belonging to the union of the boundaries, a coloring exists, since it can be rotated so that all the points become interior. Then, by the de Bruijn--Erd{\H o}s theorem, a coloring also exists for the union of the boundaries.

\section{Stability for continuous kernels}\label{sec:transfer}
Fix $n \geq 3$ and $-1<a<1$ and recall that $S$ is a unit sphere in $\R^n$.
The neighbors of $y$ in the unperturbed graph $G_a$ form the sphere
\[
N(y)=\{x\in S:\ip{x}{y}=a\}.
\]
We now show how the inequality $\chi(G_a)>n$ survives a small change in the
kernel $K$.

\begin{proof}[Proof of Theorem~\ref{thm:transfer}]
Suppose to the contrary that continuous symmetric kernels $\widetilde K_j$
converge uniformly to $K_a$ and that each $G(\widetilde K_j)$ admits a proper
$n$-coloring. Splitting classes if necessary, we may assume that all color
classes $C_{j,1},\ldots,C_{j,n}$ are nonempty. We take their closures to
obtain compact sets. These closures may contain edges, so in the argument
below we return to the color classes when choosing vertices.

The nonempty compact subsets of $S$ form a compact space under the
Hausdorff metric~\cite{Hausdorff1914,Vietoris1922}. Passing to a common subsequence, 
assume that $\overline{C_{j,i}}$ converges to $F_i$ for each $i$. 
Clearly, the sets $F_i$ form a covering of $S$. Define the \textit{palette} $P(y)$ at $y$ and the open set $O$ by
\[
 P(y)=\{i:y\in F_i\},\qquad O=S\setminus\bigcup_{i=1}^n\partial F_i .
\]
The set $O$ is open and dense, and every point of $O$ has a neighborhood on which $P$ is constant. 

\medskip\noindent\textit{A missing color on the neighbor sphere.}
For every $y\in S$ with $|P(y)| \leq n-1$, there is an $i\in P(y)$ such that
\begin{equation}\label{eq:missing}
 i\notin P(z)\qquad\text{for every }z\in N(y)\cap O .
\end{equation}
Suppose the contrary. For each $i\in P(y)$ choose a point $z_i \in N(y)\cap O$ having $i$ in its palette. 
For each color $i$, the set of admissible $z_i$ contains a nonempty subset $Z_i$ open in $N(y)$.
Each point of $N(y)$ has the form
\[
 z=ay+\sqrt{1-a^2}\,w,\qquad w\in y^\perp,\quad \|w\|=1.
\]
A generic choice of $z_i$ from the open set $Z_i$ makes the associated $w_i$ linearly independent, because $\dim y^\perp=n-1$.

Take $z_{j,i}\in C_{j,i}$ tending to $z_i$. The map $F$ defined by
\[
 Y\longmapsto\bigl(K_a(Y,z_i)\bigr)_{i\in P(y)}
\]
vanishes at $y$. The derivative $DF(y)$ is surjective because its rows on $y^\perp$ are given by the independent tangent vectors $\sqrt{1-a^2}\,w_i$. The corresponding continuous maps with entries $\widetilde K_j(Y,z_{j,i})$ converge uniformly. 
Lemma~\ref{lem:zero} applied to $F$ gives points $Y_j\to y$ satisfying $\widetilde K_j(Y_j,z_{j,i})=0$ for every $i\in P(y)$. Clearly, $Y_j \neq z_{j,i}$ for large $j$. Since every graph $G(\widetilde K_j)$ was properly colored, every $Y_j$ has color outside $P(y)$. Since there are only finitely many colors, some fixed
$\ell\notin P(y)$ satisfies $Y_j\in C_{j,\ell}$ for infinitely many $j$.
Along this subsequence, Hausdorff convergence and $Y_j\to y$
give $y\in F_\ell$, hence $\ell\in P(y)$, a contradiction. This proves \eqref{eq:missing}.

\medskip\noindent\textit{If the full palette occurs on an open set, it occurs everywhere.}
Let $W$ be the interior of $\bigcap_i F_i$ and put $N(W)=\bigcup_{x\in W}N(x)$. Clearly, $W\subset O$.
Every $y\in N(W)$ has a neighbor in $O$ with the full palette, so \eqref{eq:missing} gives $P(y)=\{1,\ldots,n\}$. The set $N(W)$ is open. Indeed, for $y\in N(x)$ a small rotation carrying $y$ to a nearby point carries $x$ to a point still in $W$.
Thus $N(W)$ is an open subset of $\bigcap_i F_i$, and hence $N(W)\subset W$.

The graph $G_a$ is connected. Two distinct unit vectors $x,z$ sufficiently close have a common neighbor. If their angular distance is $s$, the least-norm solution of $\ip{x}{u}=\ip{z}{u}=a$ has norm $|a|/\cos(s/2)<1$. Since $n \geq 3$, we can add a perpendicular vector to obtain a unit solution. Thus every graph component is open, and connectivity of the sphere gives a single component. The inclusion $N(W)\subset W$ implies that either $W=\varnothing$ or $W=S$.

\medskip\noindent\textit{The case $W=\varnothing$.}
Every $y\in O$ has $|P(y)| \leq n-1$ since its palette is locally constant. 
Consider an auxiliary coloring $c$ of $O$, in which each point $y\in O$ is assigned any color $i\in P(y)$ satisfying \eqref{eq:missing}.
If $z\in N(y)\cap O$, the color $c(y)$ does not belong to $P(z)$.
Thus $c$ is a proper coloring of the subgraph of $G_a$ induced by $O$.

By Theorem~\ref{thm:main} and the de Bruijn--Erd\H{o}s theorem~\cite[Theorem~1]{deBruijnErdos1951}, some finite subgraph $H$ of $G_a$ has chromatic number greater than $n$. There exists a rotation that places all its vertices in $O$. Indeed, for each vertex $v$, the rotations $R$ satisfying $Rv\in O$ form an open dense subset of the special orthogonal group $\mathrm{SO}(n)$, since evaluation at $v$ is an open surjection. Finitely many such sets have nonempty intersection. This contradicts the coloring of $O$.

\medskip\noindent\textit{The case $W=S$.}
Here $F_i=S$ for every color $i$. Take the model cycle and pins from
Lemma~\ref{lem:pinnedcycle}. Choose $x_{j,s,i}\in C_{j,i}$ tending to
the model pin $x_{s,i}$. On the Cartesian product of the cycle-position
spheres, impose
\[
\widetilde K_j(Y_s,x_{j,s,i}) = 0 \quad (1 \leq i \leq n-2),\qquad
\widetilde K_j(Y_s,Y_{s+1}) = 0 \quad (s\bmod m).
\]
The left-hand sides converge uniformly to the model system, whose derivative with respect to the cycle positions is invertible.
By Lemma~\ref{lem:zero}, for all sufficiently large $j$ there
are solutions $(Y_{j,s})_s$ with $Y_{j,s}\to y_s$ for every $s$.
Since the endpoints of each model edge are distinct, the same
holds for the corresponding perturbed edge for all large $j$. 
Every cycle vertex is forbidden the first $n-2$ colors. The remaining
two colors cannot color an odd cycle. This excludes $W=S$ and completes
the proof.
\end{proof}

\section{A regular finite witness}\label{sec:regularwitness}

For a finite simple graph $H$, let $L_H$ be the edge inner-product map
of Theorem~\ref{thm:stableintro}. We call a spherical realization $X$
\emph{regular} if $DL_H(X)$ has full row rank. Regularity allows us to vary
the target inner product of each edge independently.

\begin{theorem}\label{thm:finite}
Let $n \geq 3$, $-1<a<1$, and $\chi(G_a)>n$.
There is a finite $(n+1)$-chromatic graph $H$ with a regular
realization $X^0\in S^{V(H)}$ at edge inner product $a$.

Put $M_H=S^{V(H)}$ and define
\[
\Phi_H:M_H\longrightarrow\mathbb R^{E(H)},\qquad
(\Phi_H(Y))_{uv}=\ip{Y_u}{Y_v}-a.
\]
For every $\rho>0$, there is $\delta>0$ such that every continuous
map $\widetilde\Phi:M_H\to\mathbb R^{E(H)}$ satisfying
\[
\sup_{Y\in M_H}
\|\widetilde\Phi(Y)-\Phi_H(Y)\|_\infty<\delta
\]
has a zero $\widetilde Y$ with
\[
\max_{v\in V(H)}\|\widetilde Y_v-X^0_v\|<\rho.
\]
In particular, for any prescribed numbers $a_{uv}$ satisfying
$|a_{uv}-a| < \delta$, there is such a
$\widetilde Y$ with
$\ip{\widetilde Y_u}{\widetilde Y_v}=a_{uv}$ for every edge $uv$.
\end{theorem}

\begin{proof}
We construct an increasing sequence of graphs
$\mathcal H_0\subseteq\mathcal H_1\subseteq\cdots$,
assigning a position $\pi(v)\in S$ to each vertex.
The graph $\mathcal H_0$ has one isolated vertex at each point of $S$.

Given $\mathcal H_j$, construct $\mathcal H_{j+1}$ by retaining
all vertices and edges of $\mathcal H_j$ and performing every
admissible application of each of the following two types
of operations, using only vertices of $\mathcal H_j$ as pins.
A pin is a vertex whose position is held fixed while the
positions of the new vertices are chosen.
Each application introduces its own new abstract vertices,
even if their positions coincide with those of other vertices.

\begin{enumerate}
\item For every tuple of $k$ distinct old vertices,
$1 \leq k \leq n-1$, and every common neighbor $Y$ of their positions
at which the $k$ constraints have independent derivatives in $Y$,
add a vertex at $Y$ adjacent to the specified pins.
\item Use the odd integer $m$ fixed in
Lemma~\ref{lem:pinnedcycle}. For each array of old pins $p_{s,i}$,
$s\in\mathbb Z_m$, $1 \leq i \leq n-2$, consider the system
\[
\ip{Y_s}{\pi(p_{s,i})}=a,\qquad
\ip{Y_s}{Y_{s+1}}=a.
\]
For every solution $(Y_s)\in S^{\times m}$ at which the differential of this
system with respect to the cycle variables is invertible, add $m$
new vertices at the positions $Y_s$, together with the cycle edges
and the specified pin edges. Pins for distinct cycle vertices may
repeat.
\end{enumerate}

The new vertices introduced at stage $j+1$ are not used as pins
until subsequent stages. Thus applications within the same
stage do not depend on one another.

Continue this construction for all $j \geq 0$ and define
$\mathcal H$ by
\[
V(\mathcal H)=\bigcup_{j \geq 0}V(\mathcal H_j),
\qquad
E(\mathcal H)=\bigcup_{j \geq 0}E(\mathcal H_j).
\]
Every finite tuple of vertices belongs to some $\mathcal H_j$.
Consequently, every admissible application of either operation
using that tuple as pins is performed when constructing
$\mathcal H_{j+1}$.

Suppose $\mathcal H$ has a proper coloring with $n$ colors. Let $F_i$
be the closure of the set of positions of vertices of color $i$.
These closed sets cover $S$ because the initial vertices cover $S$.
Define $P$, $O$, and $W$ exactly as in the proof of
Theorem~\ref{thm:transfer}. Suppose that \eqref{eq:missing} fails at
a point $y$ with $|P(y)| \leq n-1$. As in the proof of
Theorem~\ref{thm:transfer}, choose $z_i\in N(y)\cap O$,
$i\in P(y)$, such that $i\in P(z_i)$ and the vectors
\[
w_i=\frac{z_i-ay}{\sqrt{1-a^2}}
\]
are linearly independent. 
For $u=(u_i)_{i\in P(y)}\in\mathbb R^{P(y)}$, define
\[
Y(u)=
\frac{y+\sum_{i\in P(y)}u_iw_i}
 {\left\|y+\sum_{i\in P(y)}u_iw_i\right\|}.
\]
The Jacobian with respect to $u$ of
\[
\bigl(\ip{Y(u)}{q_i}-a\bigr)_{i\in P(y)}
\]
at $u=0$ and $q_i=z_i$ is
$\sqrt{1-a^2}\,(\ip{w_h}{w_i})_{i,h\in P(y)}$, which is
invertible. The implicit function theorem therefore gives a smooth
solution $u=u(q)$ for every tuple $q=(q_i)\in S^{P(y)}$
sufficiently close to $(z_i)$, with $u((z_i))=0$ and
$u(q)\to0$ as $q_i\to z_i$. After shrinking the parameter
neighborhood if necessary, the Jacobian with respect to $u$ remains
invertible at $(u(q),q)$. Consequently, the differentials in $Y$ of
the neighbor constraints are linearly independent at $Y(u(q))$.

Choose an open neighborhood $U$ of $y$ disjoint from every $F_j$
with $j\notin P(y)$. Choose vertices $p_i$ of color $i$ whose
positions $q_i=\pi(p_i)$ are sufficiently close to $z_i$ that the
solution $Y(u(q))$ lies in $U$. Once all these pins are present,
Operation~1 adds a vertex at $Y(u(q))$ adjacent to every $p_i$.
Its color cannot belong to $P(y)$ because of these edges. Its color
cannot lie outside $P(y)$ because its position belongs to $U$.
This is a contradiction. This proves \eqref{eq:missing} for the
present coloring.

As in the proof of Theorem~\ref{thm:transfer},
\eqref{eq:missing} implies that $W$ is either empty or all of $S$.
If $W$ were empty, the missing colors would properly color $O$,
which contains a rotated copy of a finite subgraph of $G_a$ of
chromatic number greater than $n$. Thus $W=S$ and $F_i=S$ for
every $i$.

Write $(y_s)$ and $(x_{s,i})$ for the model cycle positions and
pins from Lemma~\ref{lem:pinnedcycle}. By the implicit function
theorem, there is a neighborhood
$\mathcal U\subset S^{\times m(n-2)}$ of the pin tuple $(x_{s,i})$ and
a smooth map
\[
Q=(q_{s,i})\longmapsto Y(Q)=(Y_s(Q))\in S^{\times m}
\]
such that $Y((x_{s,i}))=(y_s)$ and, for every $Q\in\mathcal U$,
\[
\ip{Y_s(Q)}{q_{s,i}}=a,\qquad
\ip{Y_s(Q)}{Y_{s+1}(Q)}=a.
\]
After shrinking $\mathcal U$, the differential of this system with
respect to the cycle variables remains invertible at $Y(Q)$ for
every $Q\in\mathcal U$.

Since $F_i=S$, choose vertices $p_{s,i}$ of color $i$ whose
positions form a tuple $Q=(\pi(p_{s,i}))_{s,i}\in\mathcal U$.
Once all these pins are present, Operation~2 adds the corresponding
odd cycle. Each cycle vertex is adjacent to pins of all colors
$1,\ldots,n-2$, so the cycle would have to be colored with the
remaining two colors. This contradiction gives
$\chi(\mathcal H)>n$.

By the de Bruijn--Erd\H{o}s theorem
\cite[Theorem~1]{deBruijnErdos1951}, $\mathcal H$ contains a finite
subgraph of chromatic number greater than $n$. Choose a vertex-minimal
such subgraph $H$. Deleting one vertex leaves an $n$-colorable graph,
so $\chi(H)=n+1$. Put $X^0=(\pi(v))_{v\in V(H)}$.
We show that $D\Phi_H(X^0)$ is surjective.

Temporarily enlarge $H$ by adjoining the complete operation that
created each of its vertices born after stage zero, including all
new vertices, edges, and pins of that operation. Repeat this for
every newly included pin born after stage zero. No operation is
associated with a stage-zero vertex. The resulting graph
$\widehat H$ is finite: every operation uses finitely many pins,
and the stage number strictly decreases whenever we pass from a
created vertex to one of its pins. Thus the depth of this recursion
is bounded by the largest stage number of a vertex of $H$.
The included operations contain every edge of $H$, since each edge
was added when at least one of its endpoints was created.

Consider the differential of the edge constraint map of
$\widehat H$ at its realization
$(\pi(v))_{v\in V(\widehat H)}$. Order the operations by stage,
choosing any order within each stage, and group their equations and
new vertex variables accordingly. Suppose there is a nontrivial
linear relation among the derivative rows. Choose the last operation
having an edge with a nonzero coefficient in this relation. On the
variables introduced by this operation, all earlier rows vanish,
while all later rows have zero coefficients in the relation.
The rows of the chosen operation are linearly independent in its
new vertex variables by construction, a contradiction.

Therefore the rows indexed by $E(H)$ are linearly independent.
These rows vanish on the variables outside $V(H)$, so they remain
independent after those variables are removed. Hence
$D\Phi_H(X^0)=DL_H(X^0)$ has full row rank.

Since $\Phi_H(X^0)=0$ and $D\Phi_H(X^0)$ is surjective,
Lemma~\ref{lem:zero} gives the asserted perturbation
property for every $\rho>0$. For prescribed edge inner products
$a_{uv}$, apply this property to
\[
(\widetilde\Phi(Y))_{uv}=\ip{Y_u}{Y_v}-a_{uv}.
\]
Indeed,
\[
\sup_{Y\in M_H}
\|\widetilde\Phi(Y)-\Phi_H(Y)\|_\infty
=\max_{uv\in E(H)}|a_{uv}-a|,
\]
finishing the proof.
\end{proof}

\begin{theorem}\label{thm:stableintro}
Fix $n \geq 3$ and $-1<a<1$. There is a finite $(n+1)$-chromatic
graph $H$ with a realization $X\in S^{V(H)}$ at edge inner
product $a$ such that
\[
L_H:S^{V(H)}\longrightarrow\mathbb R^{E(H)},\qquad
L_H(Y)=\bigl(\ip{Y_u}{Y_v}\bigr)_{uv\in E(H)}
\]
has surjective differential at $X$.
Consequently, there is an $\varepsilon>0$ such that $H$ admits a homomorphism to $G(K)$ for
every continuous symmetric kernel $K$ with
$\|K-K_a\|_\infty<\varepsilon$. In particular, $\chi(G(K))\ge n+1$.
\end{theorem}

\begin{proof} 
Put $r=(2(1-a))^{-1/2}>1/2$.
Theorem~\ref{thm:main}, after rescaling, gives $\chi(G_a)>n$.
Theorem~\ref{thm:finite} supplies a finite $(n+1)$-chromatic graph
$H$, a regular realization $X$, and the stated persistence property.

Choose $\rho>0$ such that $2\rho<\sqrt{2-2a}$, and let
$\varepsilon>0$ be the perturbation threshold supplied by
Theorem~\ref{thm:finite} for this $\rho$. For a continuous
symmetric kernel $K$, define
\[
(\Phi_K(Y))_{uv}=K(Y_u,Y_v),\qquad uv\in E(H).
\]
Then
\[
\sup_{Y\in S^{V(H)}} \|\Phi_K(Y)-\Phi_H(Y)\|_\infty \leq \|K-K_a\|_\infty.
\]
If $\|K-K_a\|_\infty<\varepsilon$, there is a zero
$\widetilde Y$ of $\Phi_K$ with
$\max_v\|\widetilde Y_v-X_v\|<\rho$. For every edge $uv$,
\[
\|\widetilde Y_u-\widetilde Y_v\|
 \geq \|X_u-X_v\|-2\rho
=\sqrt{2-2a}-2\rho>0.
\]
Thus $v\mapsto\widetilde Y_v$ is a graph homomorphism from $H$
to $G(K)$, and $\chi(G(K)) \geq \chi(H)=n+1$.

Let $f:S\to\mathbb R^n$ be continuous, put
$a=1-1/(2r^2)$, and define
\[
K_f(x,y)=\frac{1-\|f(x)-f(y)\|^2}{2r^2}.
\]
As $f$ converges uniformly to $x\mapsto rx$, the kernel $K_f$
converges uniformly to $K_a$. Since $K_f(x,y)=0$ exactly when
$\|f(x)-f(y)\|=1$, composing the homomorphism $H\to G(K_f)$
with $f$ gives a unit-distance realization of $H$ on $f(S)$.
\end{proof}

\begin{remark}[A finite witness for $1/2<r<1/\sqrt2$]\label{rem:explicit}
In each pin neighborhood $U_{s,j}$ from Lemma~\ref{lem:pinnedcycle},
apply Lemma~\ref{lem:cap} to obtain $F_{s,j}$ with distinguished set
$A_{s,j}:=V(F_{s,j})\cap U_{s,j}$. For each choice of points from $A_{s,j}$,
use Lemma~\ref{lem:pinnedcycle} to adjoin the corresponding cycle
and its pin edges to the union of the $F_{s,j}$. In a proper
$n$-coloring of this finite graph, choosing a pin of color $j$ in each
$A_{s,j}$ leaves only two colors for the cycle of the chosen tuple.
\end{remark}

\section{Discussion}
For $1/2<r<r_\Delta(n)$ the simplex cover also gives a stable upper bound
of $n+1$. At $r=r_\Delta(n)$ the diameter of each closed simplex cell
reaches one. For $n\geq 4$ it is unknown whether a sphere of larger radius can be covered by $n+1$ sets of unit diameter. Moreover, it is unclear whether for $n \geq 4$ there exists at least one value $r > r_\Delta(n)$ for which $\chi(S^{n-1}(r)) = n+1$.

The most unexpected result of the OpenAI paper on the Hadwiger--Nelson problem~\cite{OpenAI} is the transfer theorem: if an arbitrary proper $k$-coloring of $\mathbb{R}^2$ exists, then a proper $k$-coloring with measurable color classes also exists\footnote{Formally,~\cite[Theorem~1.3]{OpenAI} claims the existence of a weakly proper measurable $k$-coloring. However, a proper Lebesgue-measurable $k$-coloring can be easily derived, see~\cite[Lemma~7.1 and Corollary~8.2]{koshelev2026ordinary}.}. Denote the chromatic number with measurable color classes by $\chi_m$. At the moment, the equality $\chi(\mathbb{R}^n) = \chi_m(\mathbb{R}^n)$ was proved 
for all $n \geq 1$~\cite{OpenAI,davies2026measurable,koshelev2026ordinary}. Thus, in the case of Euclidean space, a coloring 
with everywhere dense color classes~\cite{didin2026chromatic} cannot provide a better bound than 
the measurable one. However, it remains unclear whether the analogous 
statement holds for the sphere $S^{n-1}(r)$. In particular, is it true that $\chi(S^{n-1}(r)) = \chi_m(S^{n-1}(r))$ for $n\ge3$?

\paragraph{AI use statement.} We used GPT-6.2 Astra to discover the statements that generalize our previous work and Claude Fable 5.2 to refine the exposition.
We are responsible for the mathematical content of this paper.

\paragraph{Acknowledgments.} We are grateful to Oleg Musin for introducing us to the Hopf covering theorem.  
The work of Danila Cherkashin has been accomplished with the financial support of the Bulgarian Ministry of Education and Science by the Grant No. DO1-98/26.06.2025 for NCHDC --- part of the Bulgarian National Roadmap on RIs.

\bibliographystyle{plain}
\bibliography{main}
\end{document}